\documentclass[11pt,reqno]{amsart}
\usepackage[T1]{fontenc}
\usepackage[utf8]{inputenc}
\usepackage{lmodern}
\usepackage{amsmath,amssymb,amsthm}
\usepackage{microtype}
\usepackage{xcolor}
\usepackage{hyperref}
\hypersetup{colorlinks=true,linkcolor=black,citecolor=blue!45!black,urlcolor=blue!45!black}
\newtheorem{theorem}{Theorem}
\newcommand{\R}{\mathbb{R}}
\newcommand{\E}{\mathcal{E}}
\newcommand{\proofrepo}{https://github.com/catalin-carstea/anisotropic-calderon-lean}

\title[Formal proof of anisotropic Calder\'on uniqueness]
{Formal Proof of an Anisotropic Calder\'on Uniqueness Theorem}
\author{C\u{a}t\u{a}lin I. C\^{a}rstea}
\address{Department of Applied Mathematics, National Yang Ming Chiao Tung University,
Hsinchu 30050, Taiwan}
\email{catalin.carstea@gmail.com}
\author{Lauri Oksanen}
\address{Department of Mathematics and Statistics, University of Helsinki,
P.O. Box 68, FI-00014 Helsinki, Finland}
\email{lauri.oksanen@helsinki.fi}

\begin{document}
\begin{abstract}
We announce a Lean~4 formal proof of a special case of the smooth anisotropic
Calder\'on uniqueness theorem released by OpenAI on October 6, 2026. We consider smooth
Riemannian metrics on the three-dimensional unit ball, with full boundary
measurements and Euclidean extensions outside the ball. OpenAI's release did
not include a Lean proof of their uniqueness result. Our formalization uses results from Mathlib and Tau Ceti, and it is given as a solution to a Lean comparator challenge whose check permits only the standard axioms of
Lean. The Lean proof is not a line-by-line translation of the OpenAI
manuscript.
\end{abstract}
\maketitle

\section{Introduction}

In October 2026 OpenAI released a collection of mathematical manuscripts,
several accompanied by proofs in Lean~4 \cite{openaiMath2026}. One manuscript
announces uniqueness for the smooth anisotropic Calder\'on problem from
measurements on an arbitrary nonempty open boundary patch
\cite{openaiCalderon2026}. Its theorem applies to smooth metrics on compact
manifolds of dimension at least three. The conclusion is uniqueness up to a
diffeomorphism fixing the measured patch, in accordance with the usual
coordinate invariance of the anisotropic problem
\cite{leeUhlmann1989,lassasUhlmann2001}.

We consider a more restricted statement: the manifold is the closed unit
ball in $\R^3$, measurements are available on the entire boundary, and the
metrics extend smoothly to $\R^3$, where they agree with the Euclidean
metric outside the ball.
The  main purpose of our Lean development is to verify this statement, rather than
to translate the argument of \cite{openaiCalderon2026} step by step. The autoformalization process took that manuscript as an input but was allowed to deviate from it in the proofs. 

\section{The theorem}

Consider the open unit ball
\[
  B=\{x\in\R^3:|x|<1\},
\]
and write $I$ for the Euclidean metric. An \emph{admissible metric} is a
smooth field $g:\R^3\to\R^{3\times3}$ of symmetric positive definite
matrices satisfying
\begin{equation}\label{eq:exterior}
  g(x)=I\qquad (x\in\R^3\setminus B).
\end{equation}
Thus the restriction of $g$ to $\overline B$ is a smooth Riemannian
metric whose full boundary jet agrees with that of the Euclidean metric.

We encode
boundary measurements by the minimized Dirichlet energy, which is the
quadratic form associated with the Dirichlet-to-Neumann (DN) map. The formal
theorem therefore does not require a separately defined DN map.
For a smooth function $u:\R^3\to\R$, set
\begin{equation}\label{eq:energy}
  E_g(u)
  =\int_B \nabla u(x)^{\mathsf T}\, g(x)^{-1}\nabla u(x)
     \sqrt{\det g(x)}\,dx.
\end{equation}
For a smooth function $f:\R^3\to\R$, let
\begin{equation}\label{eq:minimal}
  \E_g(f)=\inf\bigl\{E_g(u):
   u\in C^\infty(\R^3),\ u=f\text{ on }\R^3\setminus B\bigr\}.
\end{equation}
This depends only on the trace of $f$ on $\partial B$. Its value agrees with
the energy of the harmonic extension of that trace. In terms of the DN map $\Lambda_g$,
\begin{equation}\label{eq:dn}
  \E_g(f)=\bigl\langle\Lambda_g(f|_{\partial B}),
                     f|_{\partial B}\bigr\rangle_{L^2(\partial B)}.
\end{equation}
Equality of the quadratic forms $\E_g$ for a pair of metrics is equivalent, by polarization and density,
to equality of the full DN maps. 
Our formal proof \cite{carsteaOksanenFormalization2026} concerns the following theorem. 

\begin{theorem}\label{thm:main}
Let $g_1$ and $g_2$ be admissible metrics. Suppose that
\begin{equation}\label{eq:data}
  \E_{g_1}(f)=\E_{g_2}(f)
  \qquad\text{for every }f\in C^\infty(\R^3).
\end{equation}
Then there exists a smooth diffeomorphism
$\Phi:\R^3\to\R^3$, with smooth inverse, such that
\begin{equation}\label{eq:conclusion}
  \Phi(x)=x\quad \forall x\notin B, \qquad
  g_2(x)=D\Phi(x)^{\mathsf T}\,g_1(\Phi(x))\,D\Phi(x)
  \quad\forall x\in\R^3.
\end{equation}
\end{theorem}

The converse, invariance of the Dirichlet energy under such a pullback,
is the familiar diffeomorphism invariance of the anisotropic problem
\cite{leeUhlmann1989,lassasUhlmann2001}. This easy direction is formalized and proven in a comparator challenge repository~\cite{carsteaOksanenComparator2026}.

\section{Formalization}

The Lean proof \cite{carsteaOksanenFormalization2026} is given as a
solution to the comparator challenge \cite{carsteaOksanenComparator2026}.
The challenge states Theorem~\ref{thm:main} as the declaration
\[
  \texttt{OAI.AnisotropicCalderon.main}
\]
with a placeholder proof, together with the definitions it uses. The
predicate for an admissible metric is \texttt{IsAdmissibleMetric}, the
diffeomorphism condition in \eqref{eq:conclusion} is represented by
\texttt{IsIsometryFixingExterior}, and the functional in
\eqref{eq:minimal} is \texttt{minimalEnergy}. The solution proves the
same declaration in the file
\begin{center}
\small\path{OAI/Analysis/AnisotropicCalderon/Main.lean}.
\end{center}
The comparator \cite{leanComparator2026} builds the challenge and the
solution in a sandbox, checks that the solution proves the statement of
the challenge using no axioms other than \texttt{propext},
\texttt{Classical.choice} and \texttt{Quot.sound}, and replays the proof
in Lean's kernel. The challenge and its comparator configuration are kept
unchanged in the repository of the proof.

The development uses Lean~4 \cite{leanFour2021}, Mathlib \cite{mathlib2020},
and 96 source modules\footnote{The modules are distributed in the
repository with their original attribution and license notices.} from Tau Ceti \cite{tauCeti2026}. Prerequisites of the proof in \cite{openaiCalderon2026} are formalized together with its contents, even when they are well-known results. 

\section{Related work}

Other Lean formalizations of the October results include Joris
Roos's proof of an a priori $L^3\times L^3\times L^3\to L^1$ estimate
for the trilinear Hilbert transform on Schwartz functions, developed
from a blueprint based on OpenAI's argument \cite{roosTrilinear2026}.
Juan J. Manfredi formalized the main interior $C^{1,\alpha}$ estimate
for infinity-harmonic functions, together with a $C^1$ regularity
corollary and a Liouville theorem \cite{manfrediInfinity2026}.
A formalization of OpenAI's two-dimensional area-law theorem was directed
by Ruiz Chen, building in part on the formalized theory in TNLean
\cite{axiomAreaLaw2026,tnleanAreaLaw2026}.

LUOLINGc325 formalized the $11/12$ quasi-Riemann argument independently
\cite{luolingQuasiRH2026}, while Alex Kontorovich developed a direct
Lean proof of the corresponding zero-free half-plane for the Riemann
zeta function \cite{kontorovichQuasiRH2026}. Baiying Liu and LW56e2
gave Lean formalizations of refinements of OpenAI's $7/8$ zero-free
region \cite{liuZeroFree2026,lwQuasiRH2026}. Dean Cureton and Jesse
Thorner, working with Axiom Math, formalized a zero-free region for
Dedekind zeta functions \cite{axiomDedekind2026}. Peter W. Murphy used
OpenAI's formalized zero-free region in Lean proofs concerning imaginary
quadratic class numbers \cite{murphyClassNumber2026}.

Sela Navot and Jeroen Zuiddam independently developed Lean extensions
of OpenAI's matrix-multiplication exponent bound from $\mathbb C$ to
arbitrary fields \cite{navotAllFields2026,zuiddamAllFields2026}.
Jacob Sussman proved an improved power-saving bound for discrete Fourier
transforms and convolution in OpenAI's computation model
\cite{sussmanFourier2026}. Seth Douglas and Nidhal Mghirbi formalized
a non-soficity consequence of OpenAI's hyperbolic-group construction
\cite{douglasMghirbi2026}. Harald Schilly used parts of OpenAI's
formalized distance-set arguments in further Lean proofs
\cite{schillyMattila2026,schillyDistance2026}.

Separate work has reproduced released Lean proofs or compared their formal
statements with the corresponding manuscripts, including a standalone
Falconer extraction \cite{coolRmalFalconer2026}, the comparator checks and
revision records of \cite{mathcommunity2026}, and the MathVet statement-fidelity
review \cite{virk2026mathvet}. Bastounis, Circelli and Hansen analyze the
relation between natural-language claims and formal statements in connection
with an earlier Navier--Stokes announcement \cite{bastounis2026}. An
independently checked Lean argument has also related the Exponential Time
Hypothesis formalizations released by OpenAI and Anthropic
\cite{formalizeNow2026}.

OpenAI's published Lean proof in the Calder\'on family concerns
nonuniqueness for bounded measurable scalar conductivities, not the smooth
anisotropic uniqueness theorem \cite{openaiCalderonScope2026}.

\section{Checking the proof}

With Lean installed through elan, run in a clone of \url{\proofrepo}
\begin{verbatim}
lake exe cache get
lake comparator --config ComparatorChallenges/AnisotropicCalderon.json
\end{verbatim}
The comparator, which ships with Lake since Lean~4.35.0
\cite{leanRelease4350}, needs bubblewrap \cite{bubblewrap} (\texttt{bwrap})
on the path.

The challenge was set at commit \texttt{d1205fb} of
\cite{carsteaOksanenComparator2026}, on October 8, 2026, before the
autoformalization process started. To check that the challenge and the
Lean toolchain in the clone coincide with those at that commit, run
\begin{verbatim}
bash check-challenge.sh
\end{verbatim}

\section*{Acknowledgements}

We used OpenAI Codex, and to lesser extent Claude Code, in developing the Lean proof. We thank the contributors
to Lean, Mathlib and Tau Ceti for the formal libraries on which this work
relies.

C.I.C. was supported by the National Science and Technology Council of
Taiwan, grant 113-2115-M-A49-018-MY3. L.O. was supported by the European
Research Council of the European Union, grant 101086697 (LoCal), and the
Research Council of Finland, grants 347715, 353096 (Centre of Excellence
of Inverse Modelling and Imaging), and 359182 (Flagship of Advanced
Mathematics for Sensing, Imaging and Modelling). Views and opinions
expressed are those of the authors only and do not necessarily reflect
those of the European Union or the other funding organizations. Neither
the European Union nor the other funding organizations can be held
responsible for them.

\bibliographystyle{amsplain}
\bibliography{calderon-references}

@misc{openaiMath2026,
  author       = {{OpenAI}},
  title        = {Math: Mathematical manuscripts and supporting proof artifacts},
  year         = {2026},
  howpublished = {GitHub repository, first released in October},
  url          = {https://github.com/openai/math},
  note         = {Updates recorded in \url{https://github.com/openai/math/blob/main/history.md}}
}

@misc{openaiCalderon2026,
  author       = {{OpenAI}},
  title        = {Smooth anisotropic uniqueness in the {Calder\'on} problem from one boundary patch},
  year         = {2026},
  howpublished = {OpenAI Math Release preprint, dated 24 September 2026},
  url          = {https://github.com/openai/math/blob/main/preprints/Smooth-Anisotropic-Uniqueness-in-the-Calderon-Problem-from-One-Boundary-Patch-September-24-2026/paper.pdf},
  note         = {\href{https://github.com/openai/math/blob/main/preprints/Smooth-Anisotropic-Uniqueness-in-the-Calderon-Problem-from-One-Boundary-Patch-September-24-2026/paper.pdf}{Repository}}
}

@misc{openaiCalderonScope2026,
  author       = {{OpenAI}},
  title        = {Joint metric and connection recovery from one boundary patch: {Lean} formalization scope},
  year         = {2026},
  howpublished = {OpenAI mathematics repository, formalization scope for family 365},
  url          = {https://github.com/openai/math/blob/main/lean/docs/365.md},
  note         = {The listed formalization is of scalar conductivity nonuniqueness; \href{https://github.com/openai/math/blob/main/lean/docs/365.md}{Repository}}
}

@article{leeUhlmann1989,
  author       = {Lee, John M. and Uhlmann, Gunther},
  title        = {Determining anisotropic real-analytic conductivities by boundary measurements},
  journal      = {Communications on Pure and Applied Mathematics},
  volume       = {42},
  number       = {8},
  year         = {1989},
  pages        = {1097--1112},
  doi          = {10.1002/cpa.3160420804}
}

@article{lassasUhlmann2001,
  author       = {Lassas, Matti and Uhlmann, Gunther},
  title        = {On determining a {Riemannian} manifold from the {Dirichlet-to-Neumann} map},
  journal      = {Annales scientifiques de l'{\'E}cole Normale Sup{\'e}rieure},
  series       = {4},
  volume       = {34},
  number       = {5},
  year         = {2001},
  pages        = {771--787},
  doi          = {10.1016/S0012-9593(01)01076-X}
}

@misc{carsteaOksanenComparator2026,
  author       = {C\^{a}rstea, C\u{a}t\u{a}lin I. and Oksanen, Lauri},
  title        = {Anisotropic {Calder\'on} problem: {Lean} comparator challenge},
  year         = {2026},
  howpublished = {GitHub repository},
  url          = {https://github.com/l-oksanen/openai-calderon},
  note         = {The selected specification is at commit d1205fb493b3; the original main theorem has a placeholder proof; \href{https://github.com/l-oksanen/openai-calderon}{Repository}}
}

@misc{carsteaOksanenFormalization2026,
  author       = {C\^{a}rstea, C\u{a}t\u{a}lin I. and Oksanen, Lauri},
  title        = {Anisotropic {Calder\'on} uniqueness in {Lean}},
  year         = {2026},
  howpublished = {Lean source code, GitHub repository},
  url          = {https://github.com/catalin-carstea/anisotropic-calderon-lean},
  note         = {Public version examined at commit 9ed3291f0d89; \href{https://github.com/catalin-carstea/anisotropic-calderon-lean}{Repository}}
}

@misc{tauCeti2026,
  author       = {{Tau Ceti contributors}},
  title        = {Tau Ceti: A {Lean} library downstream of {Mathlib}},
  year         = {2026},
  howpublished = {Open-source formal mathematics library},
  url          = {https://github.com/TauCetiProject/TauCeti},
  note         = {The 96 files used here are bundled from revision d87e7f167605; \href{https://github.com/TauCetiProject/TauCeti}{Repository}}
}

@misc{leanComparator2026,
  author       = {{The Lean Developers}},
  title        = {Validating a {Lean} proof},
  year         = {2026},
  howpublished = {The Lean Language Reference},
  url          = {https://lean-lang.org/doc/reference/latest/ValidatingProofs/#validating-comparator},
  note         = {Section \emph{Gold standard: comparator and external checkers}, \url{https://lean-lang.org/doc/reference/latest/ValidatingProofs/#validating-comparator}}
}

@misc{leanRelease4350,
  author       = {{The Lean Developers}},
  title        = {{Lean} 4.35.0 release notes},
  year         = {2026},
  howpublished = {The Lean Language Reference},
  note         = {Released 2026-10-12; see the highlight \emph{Checking projects with lake check and lake comparator}, \url{https://lean-lang.org/doc/reference/4.36.0-rc1/releases/v4.35.0/#The-Lean-Language-Reference--Release-Notes--Lean-4___35___0-_LPAR_2026-10-12_RPAR_--Highlights--Checking-Projects-with--lake-check--and--lake-comparator}}
}

@misc{bubblewrap,
  author       = {{The Containers Organization}},
  title        = {Bubblewrap: Low-level unprivileged sandboxing tool},
  howpublished = {GitHub repository},
  url          = {https://github.com/containers/bubblewrap},
  note         = {\href{https://github.com/containers/bubblewrap}{Repository}}
}

@inproceedings{leanFour2021,
  author       = {de Moura, Leonardo and Ullrich, Sebastian},
  title        = {The {Lean} 4 theorem prover and programming language},
  booktitle    = {Automated Deduction---CADE 28},
  series       = {Lecture Notes in Computer Science},
  publisher    = {Springer},
  year         = {2021},
  doi          = {10.1007/978-3-030-79876-5_37}
}

@inproceedings{mathlib2020,
  author       = {{The mathlib Community}},
  title        = {The {Lean} Mathematical Library},
  booktitle    = {Proceedings of the 9th {ACM} {SIGPLAN} International Conference on Certified Programs and Proofs},
  series       = {CPP 2020},
  publisher    = {ACM},
  year         = {2020},
  doi          = {10.1145/3372885.3373824}
}

@misc{axiomAreaLaw2026,
  author       = {Chen, Ruiz and {AxiomMath contributors}},
  title        = {A two-dimensional area law from a global spectral gap, in {Lean} 4},
  year         = {2026},
  howpublished = {GitHub repository, released October 2026},
  url          = {https://github.com/AxiomMath/AreaLaw2D},
  note         = {Formalization of OpenAI Math Release family 265; \href{https://github.com/AxiomMath/AreaLaw2D}{Repository}}
}

@misc{tnleanAreaLaw2026,
  author       = {{TNLean contributors}},
  title        = {Tracking: two-dimensional area law from a global spectral gap},
  year         = {2026},
  howpublished = {GitHub issue and linked Lean development in TNLean},
  url          = {https://github.com/LionSR/TNLean/issues/8735},
  note         = {Reports the formalization of the area-law theorem and its corollary; \href{https://github.com/LionSR/TNLean/issues/8735}{Repository}}
}

@misc{luolingQuasiRH2026,
  author       = {{LUOLINGc325}},
  title        = {quasi-rh-lean: An independent {Lean} formalization of the $11/12$ quasi-{Riemann} hypothesis argument},
  year         = {2026},
  howpublished = {GitHub repository, first public release 28 October 2026},
  url          = {https://github.com/LUOLINGc325/quasi-rh-lean},
  note         = {\href{https://github.com/LUOLINGc325/quasi-rh-lean}{Repository}}
}

@misc{liuZeroFree2026,
  author       = {Liu, Baiying},
  title        = {Slightly improved zero-free half-planes for the quasi-{Riemann} hypothesis},
  year         = {2026},
  eprint       = {2610.12234},
  archivePrefix= {arXiv},
  primaryClass = {math.NT},
  howpublished = {arXiv:2610.12234, with accompanying Lean formalizations},
  url          = {https://arxiv.org/abs/2610.12234},
  note         = {Code: \href{https://github.com/liubaiying101/Slightly-improved-zero-free-half-planes-for-the-quasi-Riemann-hypothesis}{Lean sources}}
}

@misc{navotAllFields2026,
  author       = {Navot, Sela},
  title        = {Matrix multiplication exponent at most $9/4$ over arbitrary fields: A {Lean} formalization},
  year         = {2026},
  howpublished = {GitHub repository and Zenodo software release},
  doi          = {10.5281/zenodo.23219128},
  url          = {https://github.com/selanavot/matrix-multiplication-all-fields},
  note         = {\href{https://github.com/selanavot/matrix-multiplication-all-fields}{Repository}}
}

@misc{zuiddamAllFields2026,
  author       = {Zuiddam, Jeroen},
  title        = {$\omega(F)\leq 9/4$ over every field, in {Lean}},
  year         = {2026},
  howpublished = {GitHub repository},
  url          = {https://github.com/jzuiddam/omega-nine-quarters-all-fields},
  note         = {Independent extension of OpenAI's formalized complex-field bound; \href{https://github.com/jzuiddam/omega-nine-quarters-all-fields}{Repository}}
}

@misc{douglasMghirbi2026,
  author       = {Douglas, Seth and Mghirbi, Nidhal},
  title        = {Non-soficity of {OpenAI}'s torsion-free hyperbolic group},
  year         = {2026},
  howpublished = {Note and accompanying Lean formalization},
  doi          = {10.5281/zenodo.23202474},
  url          = {https://github.com/Apsiape/hyperbolic-group-not-sofic},
  note         = {\href{https://github.com/Apsiape/hyperbolic-group-not-sofic}{Repository}}
}

@misc{schillyMattila2026,
  author       = {Schilly, Harald},
  title        = {The {Mattila--Sj\"olin} theorem in {Lean} 4},
  year         = {2026},
  howpublished = {GitHub repository},
  url          = {https://github.com/haraldschilly/math-mattila-sjolin},
  note         = {Uses analytic lemmas from the OpenAI Falconer formalization; \href{https://github.com/haraldschilly/math-mattila-sjolin}{Repository}}
}

@misc{schillyDistance2026,
  author       = {Schilly, Harald},
  title        = {{Erd\H{o}s} similarity for distance sets: Finite patterns, a pinned obstruction, and an open window},
  year         = {2026},
  howpublished = {Preprint and Lean formalization},
  doi          = {10.5281/zenodo.23219651},
  url          = {https://github.com/haraldschilly/math-distance-set-similarity},
  note         = {\href{https://github.com/haraldschilly/math-distance-set-similarity}{Repository}}
}

@misc{coolRmalFalconer2026,
  author       = {{CoolRmal}},
  title        = {Falconer distance conjecture in all dimensions: Standalone {Lean} extraction and comparator reproduction},
  year         = {2026},
  howpublished = {GitHub repository},
  url          = {https://github.com/CoolRmal/falconer-all-dimensions},
  note         = {\href{https://github.com/CoolRmal/falconer-all-dimensions}{Repository}}
}

@misc{mathcommunity2026,
  author       = {{mathcommunity contributors}},
  title        = {oai-math-revised: Independent manuscript revisions and {Lean} comparator checks},
  year         = {2026},
  howpublished = {Collaborative GitHub project},
  url          = {https://github.com/mathcommunity/oai-math-revised},
  note         = {See \url{https://github.com/mathcommunity/oai-math-revised/blob/main/docs/LEAN_CHECKS.md}}
}

@misc{virk2026mathvet,
  author       = {Virk, Tejpal},
  title        = {{MathVet}: Independent statement-fidelity review of the {OpenAI} mathematics release},
  year         = {2026},
  howpublished = {Pre-referee review, GitHub repository and project website},
  url          = {https://math.vet/},
  note         = {See \url{https://github.com/mathvet/mathvet}}
}
\end{document}